\def\titlep{Indefinite-metric quantum field theory and operator algebra}
\documentclass[11pt]{article}
\usepackage{graphicx,ifthen}
\usepackage{amssymb}
\usepackage{amsmath}

\font\germ=eufm10 at12pt

\def\goth#1{\hbox{\germ#1}}

\newcommand{\sdag}{\scriptsize \dag}
\newcommand{\qed}{\hbox{\rule[-2pt]{3pt}{6pt}}}
\newcommand{\qedh}{\hfill\qed \\}

\newcommand{\vv}{\vspace{.3in}}

\newtheorem{Thm}{Theorem}[section]

\newtheorem{ex}[Thm]{Example}
\newtheorem{defi}[Thm]{Definition}

\newcommand{\ww}{\vv\noindent}

\def\cal#1{\mathcal #1}

\def\ltn{l_{2}({\bf N})}

\addtocounter{footnote}{1}
\def\sftt#1{
\setcounter{equation}{0}
\addtocounter{footnote}{1}
\section{#1}
}

\def\ssft#1{\subsection{#1}}

\begin{document}
%
% Personal data
%
\def\autherp{Katsunori Kawamura}
\def\emailp{e-mail: kawamura@kurims.kyoto-u.ac.jp.}
\def\addressp{College of Science and Engineering Ritsumeikan University,\\
1-1-1 Noji Higashi, Kusatsu, Shiga 525-8577,Japan
}

\newcommand{\mline}{\noindent
\thicklines
\setlength{\unitlength}{.1mm}
\begin{picture}(1000,5)
\put(0,0){\line(1,0){1250}}
\end{picture}
\par
 }
\def\sd#1{#1^{\sdag}}
\def\pca{pseudo-Cuntz algebra}
%
%%%%%%%%% Cut from here %%%%%%%%%%
%\input comm.txt
%%%%%%%%% End of Cut %%%%%%%%%
%
%
\setcounter{section}{0}
\setcounter{footnote}{0}
\setcounter{page}{1}
\pagestyle{plain}

%
% Title
%
\begin{center}
{\Large \titlep}

\ww
\autherp
\footnote{\emailp}

\noindent
{\it \addressp}
%\quad \\
%\quad \\
\quad \\
%\fbox{{\Large \today}}

\end{center}

%\maketitle
%
% Abstract
%
\begin{abstract}
It is often inevitable to introduce an indefinite-metric space
in quantum field theory.
There is a problem to determine the metric structure
of a given representation space of field operators.
We show the systematic method to determine such indefinite-metric explicitly.
At first, we choose a new 
involution $*$ of field operators instead of the original involution
$\sdag$ such that 
there is a Hilbert space $({\cal H},\langle\cdot|\cdot\rangle)$ 
with the positive-definite metric $\langle\cdot|\cdot\rangle$
which is consistent with $*$.
Next we find another hermitian form $(\cdot|\cdot)$ on ${\cal H}$
such that $({\cal H},(\cdot|\cdot))$ is a Krein space and
$(\cdot|\cdot)$ is consistent with $\sdag$.
We apply this method to various models and show that 
our results coincide with known results.
\end{abstract}


%%%%%%%%%%%%%%%%%%%%%%%%%%%%%%%%%%%%%%%%%%%%%%%%%%%%%%%
%
% Section 1
%
\sftt{Introduction}
\label{section:first}

On a complex vector space $V$, a map $(\cdot|\cdot)$ from 
$V\times V$ to ${\bf C}$
is called a {\it hermitian form} on $V$ 
if $(\cdot|\cdot)$ is sesquilinear and $\overline{(v|w)}=(w|v)$ 
for each $v,w\in V$.
In Krein space theory and quantum field theory,
it is common to call such a hermitian form an {\it inner product} or {\it metric}.
A hermitian form $(\cdot|\cdot)$ on $V$ is {\it indefinite}
if there are $v,w\in V$ such that $(w|w)<0<(v|v)$ \cite{AI,Bognar, IKL}.
Such pair $(V,(\cdot|\cdot))$ is called an {\it indefinite-metric space}
or an {\it indefinite-inner product space}.

This paper is motivated by a simple question
{\it why an indefinite-metric space appears in quantum field theory}.
According to the preface in \cite{Bognar},
the theory of indefinite-metric space 
has two origins which are relatively independent. 
One is the quantum field theory \cite{Dirac,Pauli}
and other is the functional analysis \cite{Pontrjagin,Sobolev}.
In the functional analysis \cite{AI,Bognar,IKL},
the indefinite-metric space is given at first of theory except the study 
\cite{Sobolev}
and there is no reason why an indefinite-metric space appears except
the citation from physics 

Nakanishi explained that abnormal commutation relations bring
an indefinite-metric space as the $*$-representation space of them
 ($\S$ 3, \cite{Nakanishi}).
Assume that $(\cdot|\cdot)$ is a nondegenerate hermitian form on $V$.
The {\it abnormal commutation relations} are defined by
%
% Equation 1.1
%
\begin{equation}
aa^{*}-a^{*}a=-I, \quad aa-aa=a^{*}a^{*}-a^{*}a^{*}=0.
\label{eqn:abnormalboson}
\end{equation}
Assume that $(V,(\cdot|\cdot))$ is a unital $*$-representation of 
the abnormal commutation relations and  
there is a vector $\Omega\in V$
such that $(\Omega|\Omega)>0$ and $a\Omega=0$.
Then $(a^{*}\Omega|a^{*}\Omega)=-(\Omega|\Omega)<0$.
Hence the abnormal commutation relations 
bring an indefinite-metric representation in this case.
Because the algebra generated by (\ref{eqn:abnormalboson})
is $*$-isomorphic to that of canonical commutation relations,
the reason why an indefinite-metric space appears
may be considered as the choice of vacuum.

The {\it abnormal anti-commutation relations} are defined by
%
% Equation 1.2
%
\begin{equation}
aa^{*}+a^{*}a=-I,\quad aa+aa=a^{*}a^{*}+a^{*}a^{*}=0.
\label{eqn:abnormalfermion}
\end{equation}
Assume that $a$ and $a^{*}$ are represented on $V$ such that $a^{*}$
is the adjoint of $a$ with respect to $(\cdot|\cdot)$ and the unit is preserved.
If $v\in V$ satisfies $(v|v)>0$, then $((a+a^{*})v|(a+a^{*})v)=-(v|v)<0$.
Hence the abnormal anti-commutation relations always
bring an indefinite-metric space independently in the choice of 
representation.

In this paper,
a {\it Hilbert space} $({\cal H},\langle\cdot|\cdot\rangle)$ always means
a complex vector space with a positive-definite inner product 
$\langle\cdot|\cdot\rangle$ which is complete with respect to
the norm topology induced by $\langle\cdot|\cdot\rangle$.
Any hermitian form on a complex vector space
is linear with respect to the right part 
and conjugate linear with respect to the left part.

We generalize  (\ref{eqn:abnormalboson}) and (\ref{eqn:abnormalfermion}) 
as follows. 
%
% Definition 1.1
%
\begin{defi}
\label{defi:eta}
\begin{enumerate}
%(i)
\item
A triplet $({\cal H},\langle\cdot|\cdot\rangle,\eta)$
is called a Krein triplet if 
${\cal H}$ is a Hilbert space with the inner product $\langle\cdot|\cdot\rangle$
and $\eta$ is a selfadjoint unitary on ${\cal H}$.
%(ii)
\item
For a Krein triplet $({\cal H},\langle\cdot|\cdot\rangle,\eta)$,
a family $\{a(f),a^{\sdag}(f):f\in {\cal H}\}$ satisfies
the $\eta$-canonical commutation relations ($\eta$-CCRs)
if the following holds:
\[
\left\{
\begin{array}{c}
a(f)a^{\sdag}(g)-a^{\sdag}(g)a(f)=\langle f|\eta g\rangle I,\\
\\
a(f)a(g)-a(g)a(f)=a^{\sdag}(f)a^{\sdag}(g)-a^{\sdag}(g)a^{\sdag}(f)=0
\end{array}
\right.
\quad(f,g\in {\cal H}).
\]
%(iii)
\item
For a Krein triplet $({\cal H},\langle\cdot|\cdot\rangle,\eta)$,
a family $\{a(f),a^{\sdag}(f):f\in {\cal H}\}$ satisfies
the $\eta$-canonical anti-commutation relations ($\eta$-CARs)
if the following holds:
\[
\left\{
\begin{array}{c}
a(f)a^{\sdag}(g)+a^{\sdag}(g)a(f)=\langle f|\eta g\rangle I,\\
\\
a(f)a(g)+a(g)a(f)=a^{\sdag}(f)a^{\sdag}(g)+a^{\sdag}(g)a^{\sdag}(f)=0\\
\end{array}
\right.
\quad(f,g\in {\cal H}).\]
\end{enumerate}
In both (ii) and (iii), it is understood that
the family $\{a(f),a^{\sdag}(f):f\in {\cal H}\}$ is a subset
of a unital $*$-algebra with the involution $\sdag$.
\end{defi}
A hermitian vector space $({\cal H},(\cdot|\cdot))$ is a {\it Krein space}
if there is a decomposition ${\cal H}={\cal H}_{+}
\oplus {\cal H}_{-}$ and
$({\cal H}_{\pm},\pm(\cdot|\cdot))$ is a Hilbert space \cite{AI}.
This decomposition is called the {\it fundamental decomposition} of 
$({\cal H},(\cdot|\cdot))$.
By definition,
the new hermitian form $\langle\cdot|\cdot\rangle$ on ${\cal H}$
defined by $\langle v|w\rangle\equiv (E_{+}v|E_{+}w)-(E_{-}v|E_{-}w)$ 
for $v,w\in {\cal H}$,
is positive-definite 
where $E_{\pm}$ is the projection from ${\cal H}$ onto ${\cal H}_{\pm}$.
The operator $E_{+}-E_{-}$ is called the {\it fundamental symmetry} 
of $({\cal H},(\cdot|\cdot))$.
For a Krein triplet $({\cal H},\langle\cdot|\cdot\rangle,\eta)$,
let ${\cal H}_{\pm}\equiv \{v\in {\cal H}:\eta v=\pm v\}$.
Then ${\cal H}={\cal H}_{+}\oplus {\cal H}_{-}$.
Hence $({\cal H},(\cdot|\cdot))$ is a Krein space
with the nondegenerate hermitian form $(\cdot|\cdot)$
defined by $(v|w)\equiv \langle v|\eta w\rangle$ for $v,w\in {\cal H}$.
For any operator $A$ on ${\cal H}$,
there exists unique operator $A^{\star}$ on ${\cal H}$
such that $(A^{\star}v|w)=(v|Aw)$ for each $v,w\in {\cal H}$
because $(\cdot|\cdot)$ is nondegenerate.
When $\eta=I$, the $\eta$-CCRs and the $\eta$-CARs 
coincide with ordinary CCRs and CARs ($\S$ 5.2.1 in \cite{BraRobi}), 
respectively.

%
% Theorem 1.2
%
\begin{Thm}
\label{Thm:second}
Let $({\cal H},\langle\cdot|\cdot\rangle,\eta)$ be a Krein triplet and
let ${\cal F}_{+}({\cal H})$ and ${\cal F}_{-}({\cal H})$ be the
completely symmetric and completely anti-symmetric Fock space
of the Hilbert space $({\cal H},\langle\cdot|\cdot\rangle)$, respectively.
We denote their vacuum vectors by a same symbol $\Omega$.
\begin{enumerate}
%(i)
\item
There is a Krein triplet $({\cal F}_{+}({\cal H}),\langle\cdot|\cdot\rangle,\Gamma(\eta))$
and a family $\{a(f),a^{\sdag}(f):f\in {\cal H}\}$ of operators
on ${\cal F}_{+}({\cal H})$ with an invariant dense domain ${\cal D}$
such that $\{a(f),a^{\sdag}(f):f\in {\cal H}\}$ satisfies
the $\eta$-CCRs and
\[(a^{\sdag}(f)v|w)=(v|a(f)w)\quad(v,w\in {\cal D}),\quad
a(f)\Omega=0\]
for each $f\in {\cal H}$
where $(\cdot|\cdot)$ is the hermitian form on ${\cal F}_{+}({\cal H})$
defined by $(v|w)\equiv \langle v|\Gamma(\eta)w\rangle$ 
for $v,w\in {\cal F}_{+}({\cal H})$.
%(ii)
\item
In (i), define ${\cal A}_{\eta,+}$ the $*$-algebra
generated by $a(f)$ and $a^{\sdag}(f)$ for all $f\in {\cal H}$.
Then ${\cal A}_{\eta,+}\Omega={\cal D}$.
%(iii)
\item
There is a Krein triplet $({\cal F}_{-}({\cal H}),\langle\cdot|\cdot\rangle,\Gamma(\eta))$
and a family $\{a(f),a^{\sdag}(f):f\in {\cal H}\}$ of operators
on ${\cal F}_{-}({\cal H})$ 
such that $\{a(f),a^{\sdag}(f):f\in {\cal H}\}$ satisfies
the $\eta$-CARs and
\[(a^{\sdag}(f)v|w)=(v|a(f)w)\quad(v,w\in {\cal F}_{-}({\cal H})),\quad
a(f)\Omega=0\]
for each $f\in {\cal H}$
where $(\cdot|\cdot)$ is the hermitian form on ${\cal F}_{-}({\cal H})$
defined by $(v|w)\equiv \langle v|\Gamma(\eta)w\rangle$ 
for $v,w\in {\cal F}_{-}({\cal H})$.
%(iv)
\item
In (iii), define ${\cal A}_{\eta,-}$ the $*$-algebra
generated by $a(f)$ and $a^{\sdag}(f)$ for all $f\in {\cal H}$.
Then ${\cal A}_{\eta,-}\Omega$ is dense in ${\cal F}_{-}({\cal H})$.
\end{enumerate}
Here the topologies on ${\cal F}_{+}({\cal H})$ and ${\cal F}_{-}({\cal H})$
are taken as the norm topology induced by the inner product
$\langle\cdot|\cdot\rangle$.
\end{Thm}
Theorem \ref{Thm:second} shows the following:
1) The existence of unital $*$-algebras
${\cal A}_{\eta,+}$ and ${\cal A}_{\eta,-}$ 
generated by $\eta$-CCRs and that by $\eta$-CARs, respectively.
2) The vector $\Omega$ is a cyclic vector of ${\cal A}_{\eta,\pm}$.
From the assumption of $\Omega$,
representations in Theorem \ref{Thm:second} 
are corresponded to the Fock representations of $\eta$-CCRs and $\eta$-CARs.
Especially, if $\eta=I$, then we see that
${\cal A}_{+}\equiv {\cal A}_{I,+}$ and ${\cal A}_{-}\equiv {\cal A}_{I,-}$
are the (unbounded type) CCR algebra on ${\cal F}_{+}({\cal H})$ and
 (a dense subalgebra of) the CAR algebra on ${\cal F}_{-}({\cal H})$.
We show a relation between ${\cal A}_{\pm}$ and ${\cal A}_{\eta,\pm}$.
Note that ${\cal A}_{\eta,\pm}$ is a $*$-algebra with the involution $\sdag$.
%
% Theorem 1.3
%
\begin{Thm}
\label{Thm:maintwo}
In Theorem \ref{Thm:second},
if $x^{*}$ is the adjoint operator of $x$ with respect to
$\langle\cdot|\cdot\rangle$ for $x\in {\cal A}_{\eta,\pm}$, then
\[x^{\sdag}=\Gamma(\eta)\,x^{*}\,\Gamma(\eta)^{*}\quad
\mbox{for all }x\in {\cal A}_{\eta,\pm}.\]
Especially,
${\cal A}_{\pm}=\Gamma(\eta){\cal A}_{\eta,\pm}\Gamma(\eta)^{*}$.
\end{Thm}

In algebraic quantum field theory (AQFT) \cite{Haag}, operator algebras,
that is, C$^{*}$-algebras or von Neumann algebras (W$^{*}$-algebras),
are used to describe observables in theory.
From the axiom of AQFT,
the algebra of observables is always represented on a Hilbert space.
There exist neither C$^{*}$-algebra nor W$^{*}$-algebra
containing elements which satisfy abnormal CARs.
In this way, one can not treat indefinite-metric quantum field theory (IMQFT)
in conventional AQFT apparently in the past:

\noindent
%
% Picture
%
\thicklines
\setlength{\unitlength}{.1mm}
\begin{picture}(1000,80)(-600,-20)
\put(-200,0){IMQFT}
\put(130,0){AQFT}
\put(0,0){$\Longleftrightarrow$}
\put(12,-4){\scalebox{1.5}[1.5]{$\times$}}
\put(-400,0){{\it old}:}
\end{picture}

On the other hand, our method may bring a new approach from 
AQFT to IMQFT because both CCRs and CARs are standard tools in AQFT.
By a replacement of involutions of CCRs or CARs in Theorem \ref{Thm:maintwo},
various models in IMQFT are treated explicitly:

\noindent
%
% Picture
%
\thicklines
\setlength{\unitlength}{.1mm}
\begin{picture}(1000,160)(-600,-100)
\put(-200,0){IMQFT}
\put(130,0){AQFT}
\put(0,0){$\Longleftrightarrow$}
\put(-40,-50){{\tiny replacement}}
\put(-40,-80){{\tiny of involution}}
\put(-400,0){{\it new}:}
\end{picture}

\noindent
The similar idea of the replacement of involution 
has already appeared in \cite{Pauli}.
There are many problems 
(topology, positive-definite subspace, spectral analysis of selfadjoint operator)
to treat indefinite-metric spaces and operators on them.
Our technique is quite easier than them and
the reason why the indefinite-metric space
appears is clearly explained and the relation among
$*$-algebras of field operators and their representations 
are systematically shown.

In $\S$ \ref{section:second}, we prove Theorem \ref{Thm:second} and \ref{Thm:maintwo}.
In $\S$ \ref{subsection:secondfour},
we show the difference between the $\eta$-formalism and our method.
In $\S$ \ref{section:third}, we show examples of Theorem \ref{Thm:second}
in indefinite-metric quantum field theory.
Additionally, we show representations of the BRS algebra on the Krein space 
by using the similar method of the proof of Theorem \ref{Thm:second}.
In $\S$ \ref{section:fourth}, we show the formal distribution representation
of $\eta$-commutation relations.
In Appendix \ref{section:apone}, we review ordinary CCR and CAR relations.

%%%%%%%%%%%%%%%%%%%%%%%%%%%%%%%%%%%%%%%%%%%%%%%%%%%%%%%%%%%
%
% Section 2
%
%
\sftt{Proof of Theorems}\label{section:second}
%%%%%%%%%%%%%%%%%%%%%%%%%%%%%%%%%%%%%%%%%%%%%%%%%%%%%%%%%%%%%
%
% subsection 2.1
%
\ssft{Involutive representation of involutive algebra on Krein space}
\label{subsection:secondone}
The terminology of ``$*$-algebra" is not suitable
to treat two different involutions on an algebra at once.
Hence we prepare a new terminology ``involutive algebra" instead of it.
In this paper, any algebra means an algebra over ${\bf C}$.
A map $\varphi$ on ${\cal A}$ is called an {\it involution} 
on ${\cal A}$ if $\varphi$ is a conjugate linear map which satisfies 
$\varphi(xy)=\varphi(y)\varphi(x)$ for each $x,y\in {\cal A}$ 
and $\varphi^{2}=id$.
%
% Definition 2.1
%
\begin{defi}(Chap.1, $\S$ 6, 1, \cite{Bourbaki})
A pairing $({\cal A},\varphi)$ is an involutive algebra
if $\varphi$ is an involution on an algebra ${\cal A}$.
\end{defi}

\noindent
Of course, an ordinary $*$-algebra ${\cal A}$ is an involutive algebra 
$({\cal A},*)$.
For two involutive algebras $({\cal A},\varphi)$ and $({\cal B},\psi)$, 
a homomorphism $f$ from ${\cal A}$ to ${\cal B}$ is {\it involutive}
if $f\circ \varphi=\psi\circ f$.
An automorphism $\alpha$ of $({\cal A},\varphi)$
is {\it involutive} if $\varphi\circ \alpha=\alpha\circ \varphi$.

%
% Definition 2.2
%
\begin{defi}
For a Krein triplet $({\cal H},\langle\cdot|\cdot\rangle,\eta)$,
$\pi$ is an involutive representation
of an involutive algebra $({\cal A},\varphi)$ 
if $\pi$ is a representation of ${\cal A}$
on ${\cal H}$ such that $(\pi(a)^{\star}v|w)=(v|\pi(\varphi(a))w)$
for each $a\in {\cal A}$ and $v,w\in {\cal H}$
where $\star$ is the conjugation with respect to 
the hermitian form $(\cdot|\cdot)\equiv \langle\cdot|\eta(\cdot)\rangle$.
\end{defi}
This definition is possible to be stated 
by a nondegenerate hermitian vector space
instead of ``a Krein triplet" in general.

%%%%%%%%%%%%%%%%%%%%%%%%%%%%%%%%%%%%%%%%%%%%%%%%%%%%%%%%%
%
% subsection 2.2
%
\ssft{The $\eta$-Fock space arising from a Krein triplet}
\label{subsection:secondtwo}
According to $\S$ 5.2.1 in \cite{BraRobi},
we generalize the CCR and the CAR relations 
(see also Appendix \ref{section:apone}).
Let $({\cal H},\langle\cdot|\cdot\rangle,\eta)$ be a Krein triplet and 
let ${\cal F}({\cal H})$ be the (full) Fock space of 
$({\cal H},\langle\cdot|\cdot\rangle)$.
We also denote $\langle\cdot|\cdot\rangle$ the inner product 
of ${\cal F}({\cal H})$ for the simplicity of description.
Define the operator $\Gamma(\eta)$ on ${\cal F}({\cal H})$ 
by the second quantization of $\eta$.
Then $\Gamma(\eta)$ is a selfadjoint unitary.
Hence $({\cal F}({\cal H}),\langle\cdot|\cdot\rangle,\Gamma(\eta))$ 
is a Krein triplet.
Define the hermitian form $(\cdot|\cdot)$ on ${\cal F}({\cal H})$ by
\[(x|y)\equiv \langle x|\Gamma(\eta)y\rangle\quad (x,y\in {\cal F}({\cal H})).\]
Let ${\cal F}_{+}({\cal H})$ and ${\cal F}_{-}({\cal H})$
be the Bose-Fock space and the Fermi-Fock space
with respect to the Hilbert space $({\cal H},\langle\cdot|\cdot\rangle)$.
Let $P_{\pm}$ be the projection from ${\cal F}({\cal H})$
onto ${\cal F}_{\pm}({\cal H})$.
Then we can verify that
%
% Equation 2.1
%
\begin{equation}
\label{eqn:commute}
P_{\pm}\Gamma(\eta)=\Gamma(\eta)P_{\pm}.
\end{equation}
By (\ref{eqn:commute}),
$({\cal F}_{\pm}({\cal H}),\langle\cdot|\cdot\rangle,\Gamma(\eta))$ 
is also a Krein triplet.
%
% Definition 2.3
%
\begin{defi}
Krein spaces $({\cal F}_{+}({\cal H}),(\cdot|\cdot))$  and
$({\cal F}_{-}({\cal H}),(\cdot|\cdot))$ are called 
the $\eta$-Bose-Fock space and the $\eta$-Fermi-Fock space 
by $({\cal H},\langle\cdot|\cdot\rangle,\eta)$, respectively.
\end{defi}

\noindent
In consequence, the $\eta$-Fock space is given by
the replacement of hermitian form of the ordinary Fock space.

%%%%%%%%%%%%%%%%%%%%%%%%%%%%%%%%%%%%%%%%%%%%%%%%%%%%%%%%%%%%%%%%%%
%
% subsection 2.3
%
\ssft{Proof of Theorems}
\label{subsection:secondthree}
We prove Theorem \ref{Thm:second} and \ref{Thm:maintwo} at once.
We use same assumptions for symbols ${\cal H},\langle\cdot|\cdot\rangle,(\cdot|\cdot)$,
$\eta,\Gamma(\eta)$, ${\cal F}({\cal H})$, ${\cal F}_{\pm}({\cal H})$ 
in $\S$ \ref{subsection:secondtwo}.
For each $f\in {\cal H}$, let $a(f)$ and $a^{*}(f)$ be the annihilation 
and creation operators on ${\cal F}({\cal H})$
associated with $f$ with a certain dense domain 
${\cal D}\subset {\cal F}({\cal H})$
such that $\langle a^{*}(f)v|w\rangle=\langle v|a(f)w\rangle$ for each $v,w\in {\cal D}$.
By using a symbol $\sdag$, we define the new operator $a^{\sdag}(f)$ by 
\[a^{\sdag}(f)\equiv\Gamma(\eta)a^{*}(f)\Gamma(\eta)^{*}\quad (f\in {\cal H}).\]
Then we obtain that
\[a^{\sdag}(f)=a^{*}(\eta f),\quad(a^{\sdag}(f)v|w)=(v|a(f)w)
\quad (f\in {\cal H},\,v,w\in {\cal D}).\]
Hence, $a^{\sdag}(f)$ is the adjoint operator 
of $a(f)$ with respect to the hermitian form $(\cdot|\cdot)$.
Let $a_{\pm}(f)$ and $a_{\pm}^{*}(f)$ be the annihilation 
and creation operators on ${\cal F}_{\pm}({\cal H})$.
Here the domain problem of boson case
is same with the ordinary canonical commutation relations.
Define 
\[a_{\pm}^{\sdag}(f)\equiv \Gamma(\eta)a_{\pm}^{*}(f)\Gamma(\eta)^{*}
\quad(f\in {\cal H}).\]
By (\ref{eqn:commute}),
we see that $a^{\sdag}_{\pm}(f)$ is the adjoint operator 
of $a_{\pm}(f)$ with respect to $(\cdot|\cdot)$.
From $a_{\pm}(f)a^{*}_{\pm}(g)\mp a^{*}_{\pm}(g)a_{\pm}(f)=\langle f|g\rangle I$,
one computes straightforwardly that
$\{a_{+}(f),a_{+}^{\sdag}(f):f\in {\cal H}\}$
satisfies $\eta$-CCRs and 
$\{a_{-}(f),a_{-}^{\sdag}(f):f\in {\cal H}\}$ satisfies $\eta$-CARs.
The remaining statements in Theorem \ref{Thm:second} and \ref{Thm:maintwo} hold
by the properties of $\{a_{\pm}(f),a^{*}_{\pm}(f):f\in {\cal H}\}$
on ${\cal F}_{\pm}({\cal H})$.
\qedh

\noindent
By construction, both $\eta$-CCRs and $\eta$-CARs are obtained by
the replacement of involution of the ordinary CCRs and CARs, respectively.
We call their representations on $({\cal F}_{+}({\cal H}),(\cdot|\cdot))$ and
$({\cal F}_{-}({\cal H}),(\cdot|\cdot))$ by
the {\it $\eta$-Bose-Fock representation} of $\eta$-CCRs and
the {\it $\eta$-Fermi-Fock representation} of $\eta$-CARs, respectively.
We simply call them by the {\it $\eta$-Fock representations}.

%%%%%%%%%%%%%%%%%%%%%%%%%%%%%%%%%%%%%%%%%%%%%%%%%%%%%%%%%%%%%%%
%
% subsection 2.4
%
\ssft{Difference between $\eta$-formalism and our method}
\label{subsection:secondfour}
The $\eta$-formalism (unitary trick) is the old-fashioned way
of treating the indefinite-metric Hilbert space. 
In $\S$ 8 (B) of \cite{Nakanishi}, 
it is written that $\eta$-formalism depends on the choice of particular 
pseudo-orthonormal basis and  this fact is very inconvenient.
The aim of $\eta$-formalism is to replace the indefinite metric
to a positive-definite metric without any consideration 
of the involution of field operators.
After the introduction of the positive-definite metric associated with 
the metric operator $\eta$, the new involution is introduced.
Because it is necessary for the $\eta$-formalism to know the indefinite-metric
of a given theory, the $\eta$-formalism itself can not determine
the indefinite metric.
In fact, if one applies the $\eta$-formalism of 
examples in $\S$ \ref{section:third},
then one must compute the complete system of 
a given representation in order to determine $\eta$.

On the other hand, the first aim of our method is the replacement
of the involution of algebra of field operators.
We consider that the almost indefinite-metric quantum field theory 
is caused by a priori involution in theory.
Our method replaces such ill involution to the authorized involution
in the theory of operator algebra as like as the algebra ${\cal A}_{CCR}$ 
of CCRs and the algebra ${\cal A}_{CAR}$ of CARs. 
Such method is possible to be applied when
the commutation relations associated the ill involution
is the type of $\eta$-CCRs or $\eta$-CARs.
After the replacement of involution, we can apply the 
theory of the ordinary Fock representation theory of 
${\cal A}_{CCR}$ and ${\cal A}_{CAR}$.
These theories are assured by the theory of operator algebra \cite{BraRobi}.
At the last, we place back the involution and determine
the indefinite metric at once.
Our method clarifies the indefinite-metric structure of 
quantum field theory in the point of view of operator algebra 
which is based on the theory of Hilbert spaces with positive-definite
metric explicitly.
%%%%%%%%%%%%%%%%%%%%%%%%%%%%%%%%%%%%%%%%%%%%%%%%%%%%%%%%%%%
%
% Section 3
%
\sftt{Examples}\label{section:third}
We use same assumptions for symbols 
${\cal H},\langle\cdot|\cdot\rangle,(\cdot|\cdot),\eta,\Gamma(\eta)$,
${\cal F}_{\pm}({\cal H})$, $a_{\pm}(f)$, $a_{\pm}^{\sdag}(f)$ 
in $\S$ \ref{subsection:secondthree}.
%
% Example 3.1
%
\begin{ex}{\rm
Let $\eta\equiv -I$.
Then we obtain that
\[a_{\pm}(f)a_{\pm}^{\sdag}(g)\mp a^{\sdag}_{\pm}(g)a_{\pm}(f)=-\langle f|g\rangle I
\quad(f,g\in {\cal H}).\]
The $\eta$-CCR and $\eta$-CAR relations coincide with 
(\ref{eqn:abnormalboson}) and (\ref{eqn:abnormalfermion}), respectively 
when ${\rm dim}{\cal H}=1$.
Fix a completely orthonormal basis
$\{e_{n}\}_{n\in \Lambda}$ of $({\cal H},\langle\cdot|\cdot\rangle)$
and define
\[a_{n,\pm}\equiv a_{\pm}(e_{n}),\quad a^{\sdag}_{n,\pm}=a^{\sdag}_{\pm}(e_{n})
\quad(n\in \Lambda).\]
Then we obtain that
\[a_{n,\pm}a^{\sdag}_{m,\pm}\mp a^{\sdag}_{m,\pm}a_{n,\pm}=-\delta_{nm}I,\]
\[a_{n,\pm}a_{m,\pm}\mp a_{m,\pm}a_{n,\pm}=a_{n,\pm}^{\sdag}a_{m,\pm}^{\sdag}
\mp a_{m,\pm}^{\sdag}a_{n,\pm}^{\sdag}=0
\quad(n,m\in\Lambda).\]
Such $\eta$-CARs appear in the Lee model
(See \cite{Nakanishi}, $\S$ 12,
see also \cite{Fuda, Heisenberg,Kaller}) when ${\rm dim}{\cal H}=1$.
}
\end{ex}
%
% Example 3.2
%
\begin{ex}{\rm
Let $(\ltn,\langle\cdot|\cdot\rangle)$ be the Hilbert space 
with the standard basis $\{e_{n}\}_{n\in {\bf N}}$.
Define ${\cal H}\equiv {\cal H}_{+}\oplus {\cal H}_{-}$
for ${\cal H}_{\pm}\equiv \ltn$
and define the basis $\{e_{n,\pm}\}_{n\in{\bf N}}$
of ${\cal H}_{\pm}$ by $e_{n,\pm}\equiv e_{n}$ for each $n$.
\begin{enumerate}
%(i)
\item
Define $\eta$ and operators $\alpha_{n}$ and $\beta_{n}$  by
\[\eta e_{n,\pm}\equiv e_{n,\mp},\quad
\alpha_{n}\equiv a_{+}(e_{n,+}),\quad \beta_{n}\equiv a_{+}(e_{n,-})
\quad (n\in{\bf N}).\]
Then we obtain
\[
\left\{
\begin{array}{c}
\alpha_{n}\beta^{\sdag}_{m}-\beta^{\sdag}_{m}\alpha_{n}=
\beta_{n}\alpha^{\sdag}_{m}-\alpha^{\sdag}_{m}\beta_{n}=\delta_{n,m}I,\\
\\
\alpha_{n}\alpha^{\sdag}_{m}-\alpha^{\sdag}_{m}\alpha_{n}=
\beta_{n}\beta^{\sdag}_{m}-\beta^{\sdag}_{m}\beta_{n}=0
\end{array}
\right.
\quad (n,m\in{\bf N}).\]
Other commutators vanish.
These commutation relations appear in the Froissart model
($\S$ 13 in \cite{Nakanishi}, see also \cite{Froissart}).
%(ii)
\item
Define $\eta$ and operators $a_{n}$ and $a_{n}^{\sdag}$ by
\[\eta e_{n,\pm}\equiv \pm e_{n,\pm}\quad
a_{2n}\equiv a_{-}(e_{n,+}),\quad 
a_{2n-1}\equiv a_{-}(e_{n,-})\quad (n\in {\bf N}).\]
Then we obtain
\[a_{n}a_{m}^{\sdag}+a_{m}^{\sdag}a_{n}=(-1)^{n}\delta_{n,m}I,\quad
a_{n}a_{m}+a_{m}a_{n}=a_{n}^{\sdag}a_{m}^{\sdag}+a_{m}^{\sdag}a_{n}^{\sdag}=0\]
for $n,m\in {\bf N}$.
These are called ICARs in $\S$ 4 of \cite{AK05}.
By using these, we constructed FP (anti) ghosts in string theory.
\end{enumerate}
}
\end{ex}
%
% Example 3.3
%
\begin{ex}{\rm
Assume that ${\rm dim}{\cal H}=2$ and $e_{1},e_{2}$ are 
orthonormal basis of ${\cal H}$.
A selfadjoint unitary $\eta$ on ${\cal H}$ is one of the following:
\[\pm I,\quad 
\eta(\theta,\xi)\equiv 
\left(
\begin{array}{cc}
\cos\xi &e^{\sqrt{-1}\theta}\sin\xi\\
e^{-\sqrt{-1}\theta}\sin\xi&-\cos\xi\\
\end{array}
\right)\quad (\xi,\theta\in [0,2\pi)).
\]
Define operators $a_{1}\equiv a_{+}(e_{1})$ and $a_{2}\equiv a_{+}(e_{2})$
and define ${\cal A}(\eta)$ the involutive algebra generated by $a_{1},a_{2}$.
Then we see that ${\cal A}(\eta)$ is $*$-isomorphic to ${\cal A}(I)$ 
for each $\eta$.
In consequence, there exists an algebras of $\eta$-CCRs
uniquely up to isomorphism when ${\rm rank}\, \eta=2$.

On the other hand,
if ${\cal B}(\eta)$ is the $\eta$-CAR algebra, then 
${\cal B}(I)$ and ${\cal B}(-I)$ are not involutively isomorphic.
Hence we see that there exist two mutually non-isomorphic 
algebras of $\eta$-CARs at least when ${\rm rank}\, \eta=2$.
}
\end{ex}
%
% Example 3.4
%
\begin{ex}{\rm
The Minkowski metric $g$ often appears in commutation relations.
We show that such case is a special $\eta$-CCRs.
Define $g=(g_{\mu\nu})_{\mu,\nu=0,1,2,3}
={\rm diag}(1,-1,-1,-1)\in M_{4}({\bf R})$.
In the quantum electromagnetic dynamics
at the Feynman gauge $\alpha=1$ ($\S$ 2.3, \cite{NO}),
the following commutation relations appear:
%
% Equation 3.1
%
\begin{equation}
\label{eqn:feynmann}
a_{\mu}a_{\nu}^{\sdag}-a_{\nu}^{\sdag}a_{\mu}=-g_{\mu\nu}I\quad(\mu,\nu=0,1,2,3)
\end{equation}
where we omit the suffix of operators from the originals
except $\mu,\nu$ because other suffix does not 
bring abnormal commutation relations.
We reformulate these commutation relations by $\eta$-CCR as follows:
Let ${\cal H}\equiv {\bf C}^{4}
={\bf C}e_{0}\oplus {\bf C}e_{1}\oplus {\bf C}e_{2}\oplus {\bf C}e_{3}$
with the standard inner product $\langle\cdot|\cdot\rangle$ of ${\bf C}^{4}$ 
with respect to $e_{0},e_{1},e_{2},e_{3}$, and let 
\[\eta\equiv -g.\]
Define $a_{\mu}\equiv a_{-}(e_{\mu})$ for $\mu=0,1,2,3$.
Then we obtain (\ref{eqn:feynmann})
on the $\eta$-Bose-Fock space ${\cal F}_{+}({\bf C}^{4})$.
We show the fundamental decomposition of ${\cal F}_{+}({\bf C}^{4})$.
Define $V_{+}\equiv {\bf C}e_{1}\oplus {\bf C}e_{2}\oplus {\bf C}e_{3}$ and
\[V_{n,m}\equiv V_{+}^{\vee n}\vee {\bf C}e_{0}^{\vee m}\quad(n+m\geq 1)\]
where $\vee$ is the symmetric tensor product.
Define
\[{\cal F}_{+}({\bf C}^{4})_{+}\equiv 
{\bf C}\Omega\oplus \bigoplus_{n}\bigoplus_{k+l=n,\, l:even}V_{k,l},\quad
{\cal F}_{+}({\bf C}^{4})_{-}\equiv
\bigoplus_{n}\bigoplus_{k+l=n,\, l:odd}V_{k,l}.\]
We see that $({\cal F}_{+}({\bf C}^{4})_{\pm},\pm(\cdot|\cdot))$ 
is positive-definite where 
$(v|w)\equiv \langle v|\Gamma(\eta)w\rangle$ for $v,w\in {\cal F}_{+}({\bf C}^{4})$.
}
\end{ex}

\noindent
Remark that 
the eigenvalue problem in each model is not changed by the choice of metric.
Hence our manipulation does not change any eigenvalue problem.

%
% Example 3.5
%
\begin{ex}{\rm
The BRS algebra (\cite{NO} $\S$ 3.4.2)
is an involutive algebra $({\cal A},\sdag)$
generated by $Q_{B},Q_{C}$ such that
%
% Equation 
%
\begin{equation}
\label{eqn:brs}
Q_{B}^{\sdag}=Q_{B},\quad Q_{B}^{2}=0,\quad
Q_{C}^{\sdag}=Q_{C},\quad Q_{C}Q_{B}-Q_{B}Q_{C}=-\sqrt{-1}Q_{B}.
\end{equation}
It is known that a nondegenerate involutive 
representation of ${\cal A}$ is an indefinite-metric space
because $(Q_{B}v|Q_{B}v)=(v|Q_{B}^{2}v)=0$ for each vector $v$.
Let $V$ be a cyclic representation space of ${\cal A}$ 
with a cyclic vector $\Omega$.
If $Q_{B}\Omega=0$, then $Q_{B}Q_{C}\Omega=0$.
Hence this representation is degenerate.

For $a\in {\bf R}$, we define a nondegenerate representation of ${\cal A}$
on $V\equiv {\bf C}^{2}$.
Define matrices $Q_{B},Q_{C},U\in M_{2}({\bf C})$ by
\[Q_{B}\equiv 
\left(
\begin{array}{cc}
0&1\\
0&0\\
\end{array}
\right),\quad
Q_{C}\equiv 
\left(
\begin{array}{cc}
a+\sqrt{-1}/2&0\\
0&a-\sqrt{-1}/2\\
\end{array}
\right),
\quad
U\equiv 
\left(
\begin{array}{cc}
0&1\\
1&0\\
\end{array}
\right).
\]
Define the new involution $\sdag$ on $M_{2}({\bf C})$ by
\[x^{\sdag}\equiv Ux^{*}U^{*}\quad(x\in M_{2}({\bf C}))\]
where $*$ is the hermite conjugation on $M_{2}({\bf C})$.
Then we can verify that $Q_{B}$ and $Q_{C}$ satisfy (\ref{eqn:brs}).
Define the hermitian form $(\cdot|\cdot)$ on ${\bf C}^{2}$
by $(v|w)\equiv \langle v|Uw\rangle$ for $v,w\in {\bf C}^{2}$.
Then $(x^{\sdag}v|w)=(v|xw)$ for each $x\in M_{2}({\bf C})$ 
and  $v,w\in {\bf C}^{2}$.
Therefore the Krein space $({\bf C}^{2},(\cdot|\cdot))$ is a nondegenerate
involutive representation of $({\cal A},\sdag)$.
Let ${\cal H}_{\pm}\equiv 
\{v\in {\bf C}^{2}:Uv=\pm v\}$.
Then
${\cal H}_{+}={\bf C}(e_{1}+e_{2})$ and
${\cal H}_{-}={\bf C}(e_{1}-e_{2})$
where $e_{1},e_{2}$ are standard basis of ${\bf C}^{2}$.
}
\end{ex}
%%%%%%%%%%%%%%%%%%%%%%%%%%%%%%%%%%%%%%%%%%%%%%%%%%%%%%%%%%%%%
%
% Section 4
% 
\sftt{Distribution representation of $\eta$-commutation relations}
\label{section:fourth}
Let $X\equiv {\bf R}^{n}$.
In theoretical physics,
field operators are usually described by operator-valued distribution 
as like as 
\[[a(p), a^{\sdag}(q)]=\delta(p-q)I\quad (p,q\in X)\]
where $[x,y]\equiv xy-yx$.
We show $\eta$-CCRs according to the above notation formally.
Let $L_{2}(X)$ be the Hilbert space of all square-integrable complex-valued 
functions on $X$ with the $L_{2}$-inner product $\langle\cdot|\cdot\rangle$.
Let $\eta$ be a selfadjoint unitary on $L_{2}(X)$.
Assume that we can denote {\it formally} that 
\[(\eta f)(x)=\int \!dp\, f(p)\eta(p,x)\quad (x\in X)\]
for each suitable function $f$ on $X$.
Then the following holds:
\[\overline{\eta(p,q)}=\eta(q,p),\quad
\int \!dp\,\eta(q,p)\eta(p,x)=\delta(x-q).\]
If $[a(f),a^{\sdag}(g)]=\langle f|\eta g\rangle I$ for each $f,g\in L_{2}(X)$, then
\[[a(p), a^{\sdag}(q)]=\eta(p,q)I\quad (p,q\in X).\]
These are verified {\it formally} by using only algebraic manipulation.

%
% Example 4.1
%
\begin{ex}{\rm
Consider two operator-valued distributions $a(p)$ and $b(p)$ 
which satisfy the following:
\[[a(p),b^{\sdag}(q)]=[b(p),a^{\sdag}(q)]=\delta(p-q)I,\quad
[a(p),a^{\sdag}(q)]=[b(p),b^{\sdag}(q)]=0\quad(p,q\in X).\]
Let $c(p,1)\equiv a(p)$ and $c(p,2)\equiv b(p)$.
On ${\cal H}\equiv L_{2}(X)\otimes {\bf C}^{2}$,
define the correspondence ${\cal H}\ni f\otimes v\mapsto c(f\otimes v)$ by
\[c(f\otimes v)=\sum_{i=1}^{2}\int\! dp\, c(p,i)\overline{f(p)}\overline{v_{i}}
\quad (f\in L_{2}(X),\, v=(v_{1},v_{2})\in {\bf C}^{2}).\]
Hence the $\eta$-CCRs are given as
\[[c(f\otimes v),\,\,c^{\sdag}(g\otimes w)]=\langle f\otimes v|(I\otimes \eta_{0})
g\otimes w\rangle I \]
where $\eta_{0}\equiv 
\left(\begin{array}{cc}
0&1\\
1&0\\
\end{array}
\right)$.
Define $\eta\equiv I\otimes \eta_{0}$ and identify an operator on ${\cal H}$
as a $2\times 2$-matrix consisting of operators on $L_{2}(X)$.
Then we obtain a distribution type $\eta$-CCRs,
\[[c(p,i),\,\,c^{\sdag}(p,j)]=(\eta(p,q))_{i,j}I\quad(i,j=1,2)\]
where	
\[\eta(p,q)=\delta(p-q)\otimes \eta_{0}=
\left(\begin{array}{cc}
0&\delta(p-q)\\
\delta(p-q)&0\\
\end{array}
\right).\]
}
\end{ex}

\ww
{\bf Acknowledgement:}
The author would like to express my sincere thanks to Noboru Nakanishi
and Takeshi Nozawa for beneficial advices for this article.

\appendix
%%%%%%%%%%%%%%%%%%%%%%%%%%%%%%%%%%%%%%%%%%%%%%%%%%%%%%%%%%%%%%%%%
%
% Appendix A
%
\sftt{The CCR and CAR relations}
\label{section:apone}
We review the CCR and CAR relations according to
$\S$ 5.2.1 in \cite{BraRobi}.
Let ${\cal H}$ be a Hilbert space with the inner product $\langle\cdot|\cdot\rangle$
and let ${\cal H}^{\otimes n}$ denote the $n$-fold tensor product of ${\cal H}$
with itself.
Further introduce the (full) Fock space ${\cal F}({\cal H})$ by
\[{\cal F}({\cal H})\equiv {\bf C}\Omega\oplus {\cal H}
\oplus {\cal H}^{\otimes 2}\oplus \cdots\]
where $\Omega$ is a unit vector.
Define projections $P_{\pm}$ on ${\cal F}({\cal H})$ by 
$P_{\pm}\Omega\equiv \Omega$,
\[P_{+}(f_{1}\otimes \cdots\otimes f_{n})\equiv 
(n!)^{-1}\sum_{\sigma\in{\goth S}_{n}}f_{\sigma(1)}\otimes\cdots
\otimes f_{\sigma(n)},\]
\[P_{-}(f_{1}\otimes \cdots\otimes f_{n})\equiv 
(n!)^{-1}\sum_{\sigma\in{\goth S}_{n}}{\rm sgn}(\sigma)f_{\sigma(1)}
\otimes\cdots\otimes f_{\sigma(n)}\]
for all $f_{1},\ldots,f_{n}\in {\cal H}$.
Subspaces ${\cal F}_{+}({\cal H})\equiv P_{+}{\cal F}({\cal H})$ 
and ${\cal F}_{-}({\cal H})\equiv P_{-}{\cal F}({\cal H})$ 
are called the {\it Bose-Fock space} and the {\it Fermi-Fock space}, 
respectively.
For a unitary $U$ on ${\cal H}$,
let $U^{\otimes n}$ denote the $n$-fold tensor product of $U$ with itself.
The unitary operator $\Gamma(U)$ on ${\cal F}({\cal H})$ defined by
\[\Gamma(U)\equiv I_{{\bf C}\Omega}\oplus U\oplus U^{\otimes 2}\oplus \cdots\]
is called the {\it second quantization} of $U$.
For $f\in {\cal H}$,
define  operators $a(f)$ and $a^{*}(f)$ on ${\cal F}({\cal H})$ by
$a(f)\Omega\equiv 0$, $a^{*}(f)\Omega\equiv f$ and
\[a(f)(f_{1}\otimes \cdots\otimes f_{n})\equiv n^{1/2}\langle f|f_{1}\rangle
f_{2}\otimes \cdots\otimes f_{n},\]
\[a^{*}(f)(f_{1}\otimes \cdots\otimes f_{n})\equiv (n+1)^{1/2}
f\otimes f_{1}\otimes \cdots\otimes f_{n}.\]
On the domain ${\cal D}\equiv \{(v_{n})_{n\geq 0}\in {\cal F}({\cal H}):
\sum_{n\geq 0}n\|v_{n}\|^{2}<\infty\}$,
one has the adjoint relation 
$\langle a^{*}(f)v|w\rangle=\langle v|a(f)w\rangle$ for $v,w\in {\cal D}$.
Define
\[a_{\pm}(f)\equiv P_{\pm}a(f)P_{\pm},\quad
a^{*}_{\pm}(f)\equiv P_{\pm}a^{*}(f)P_{\pm}\quad(f\in {\cal H}).\]
Then we obtain that
\[a_{+}(f)a^{*}_{+}(g)- a^{*}_{+}(g)a_{+}(f)=\langle f|g\rangle I,\]
\[a_{+}(f)a_{+}(g)-a_{+}(g)a_{+}(f)=
a^{*}_{+}(f)a^{*}_{+}(g)- a^{*}_{+}(g)a^{*}_{+}(f)=0,\]
\[a_{-}(f)a^{*}_{-}(g)+ a^{*}_{-}(g)a_{-}(f)=\langle f|g\rangle I,\]
\[a_{-}(f)a_{-}(g)+ a_{-}(g)a_{-}(f)=
a^{*}_{-}(f)a^{*}_{-}(g)+ a^{*}_{-}(g)a^{*}_{-}(f)
=0\quad(f,g\in {\cal H}).\]
The equations for $a_{+}(f)$ and $a_{+}^{*}(g)$
are called the {\it canonical commutation relations (CCRs)}
and those for $a_{-}(f)$ and $a_{-}^{*}(g)$
are called the {\it canonical anti-commutation relations (CARs)}.

%%%%%%%%%%%%%%%%%%%%%%%%%

\end{document}